\documentclass[12pt]{article} 
\usepackage{index}
\usepackage[reqno]{amsmath}
\usepackage{fixmath}
\usepackage{amssymb, amsthm, enumerate}
\usepackage{mathtools}

\usepackage{hyperref}

\DeclarePairedDelimiter{\abs}{\lvert}{\rvert}
\DeclarePairedDelimiter{\norm}{\lVert}{\rVert}

\newtheoremstyle{plainsl}%
	{\topsep}
	{\topsep}
	{\slshape} 
	{}
	{\normalfont\bfseries}
	{.}
	{ }
	{}

\swapnumbers

{\theoremstyle{plainsl}
\newtheorem{theorem}{Theorem}[section]
\newtheorem{lemma}[theorem]{Lemma}
\newtheorem{corollary}[theorem]{Corollary}}
{\theoremstyle{remark}
}

\newcommand\cref[1]{Corollary~\ref{cor:#1}}

\renewcommand\proof{\noindent\textsl{Proof. }}
\newcommand\sqr[2]{{\vbox{\hrule height.#2pt
    \hbox{\vrule width.#2pt height#1pt \kern#1pt
        \vrule width.#2pt}\hrule height.#2pt}}}
\renewcommand\qed{%
	\ifmmode\eqno\sqr53
	\else\nolinebreak\ \hfill\sqr53\medbreak\fi}

\numberwithin{equation}{section}

\newcommand\al{\alpha}
\newcommand\be{\beta}
\newcommand\de{\delta}
\newcommand\De{\Delta}
\newcommand\eps{\epsilon}
\newcommand\ga{\gamma}
\newcommand\Ga{\Gamma}

\newcommand\la{\lambda}

\renewcommand\th{\theta} 

\newcommand\cx{{\mathbb C}}

\newcommand\flde{{\mathbb E}}
\newcommand\ints{{\mathbb Z}}
\newcommand\re{{\mathbb R}}
\newcommand\rats{{\mathbb Q}}
\newcommand\comp[1]{{\mkern2mu\overline{\mkern-2mu#1}}}
\newcommand\diff{\mathbin{\mkern-1.5mu\setminus\mkern-1.5mu}}

\newcommand\seq[3]{#1_{#2},\ldots,#1_{#3}}

\newcommand\pmat[1]{\begin{pmatrix} #1 \end{pmatrix}}

\DeclareMathOperator{\rk}{rk}
\DeclareMathOperator{\tr}{tr}

\DeclareMathOperator{\arcs}{arcs}
\DeclareMathOperator{\esupp}{esupp}

\newcommand\ip[2]{\langle #1,#2\rangle}

\newcommand\alg[1]{\langle#1\rangle}

\title{Real State Transfer} 

\author{
	Chris Godsil\footnote{
		Research supported by Natural Sciences and Engineering Council of Canada, 
		Grant No. RGPIN-9439}\\
	Combinatorics \& Optimization\\
	University of Waterloo}

\begin{document}
\maketitle
	
\begin{abstract}
    A continuous quantum walk on a graph $X$ with adjacency matrix $A$ is specified
	by the 1-parameter family of unitary matrices $U(t)=\exp(itA)$. These matrices act
	on the state space of a quantum system, the states of which we may represent by
	density matrices, positive semidefinite matrices with rows and columns indexed by 
	$V(X)$ and with trace $1$. The square of the absolute values of the entries of
	a column of $U(t)$ define a probability density on $V(X)$, and it is precisely these
	densities that predict the outcomes of measurements. There are two special cases of 
	physical interest: when the column density is supported on a vertex, and when it is
	uniform. In the first case we have perfect state transfer; in the second, uniform mixing.
	
	There are many results concerning state transfer and uniform mixing. In this paper we show
	that these results on perfect state transfer hold largely because at the time
	it occurs, the density matrix is real. We also show that the results on uniform mixing
	obtained so far hold because the entries of the density matrix are algebraic numbers.
	As a consequence of these we derive strong restrictions on the occurence of uniform mixing
	on bipartite graphs and on oriented graphs.
\end{abstract}

\section{Introduction}
\label{sec:intro}

Let $X$ be a graph with adjacency matrix $A$. The states of a 
\textsl{continuous quantum walk} on $X$ are represented by
positive semidefinite matrices and with trace $1$.
Physicists refer to such matrices as \textsl{density matrices}. 

Let $X$ be a graph with adjacency matrix $A$. The states of a \textsl{continuous quantum walk} 
on $X$ are density matrices with rows and columns indexed by $V(X)$. 
The behaviour of the walk is determined by its initial state and the \textsl{transition matrix} $U(t)$, 
defined by
\[
	U(t) = \exp(itA).
\]
This is a symmetric and unitary matrix. If the initial state of our walk is given by
a density matrix $D$, then the state $D(t)$ at time $t$ is given by
\[
	D(t) = U(t)DU(-t).
\]
A density matrix $D$ represents a \textsl{pure state} if $\rk(D)=1$. In this case
$D=zz^*$ for some complex unit vector $z$. If $e_a$ denotes the standard basis
vector in $\cx^{V(X)}$ indexed by the vertex $a$, then
\[
	D_a = e_ae_a^T
\]
is the pure state associated to the vertex $a$.

One question of interest to physicists is whether, for two distinct vertices $a$
and $b$, there is a time $t$ such that $D_a(t)=D_b$. In this case we say that
there is \textsl{perfect state transfer} from $a$ to $b$.
If perfect state transfer occurs, a number of interesting consequences have been extablished.
(See e.g. \cite{st-transfer}.)

We will summarize these, after developing some more terminology. We assume that the
eigenvalues of $A$ are $\seq\th1m$, and that the matrix $E_r$ represents orthogonal
projection onto the $\th_r$-eigenspace. Then $E_rE_r=\de_{r,s}E_r$ and $\sum_r E_r=I$
and we have the spectral decomposition 
\[
	A = \sum_r \th_r E_r.
\]
(One consequence of this is that $U(t)=\sum_r e^{it\th_r}E_r$.) Vertices $a$ and $b$ of $X$ are
said to be \textsl{cospectral} if the graphs $X\diff a$ and $X\diff b$ are cospectral, that is,
they have the same characteristic polynomial. It is known that $a$ and $b$ are cospectral
if and only if
\[
	\norm{E_re_a} = \norm{E_re_b}
\]
for all $r$. We say that $a$ and $b$ are \textsl{strongly cospectral} if, for $r$,
\[
	E_re_a = \pm E_re_b.
\]
For more about cospectral and strongly cospectral vertices, see \cite{cgjs-strong}.

We can now list the consequences of perfect state transfer. If there is pst from $a$ to $b$ at 
time $t$, then:
\begin{enumerate}[(a)]
	\item 
	There is pst from $b$ to $a$ at time $t$.
	\item
	$D_a(2t)=D_a$.
	\item
	If $E_kPE_\ell\ne0$ and $E_rPE_s\ne0$ and $k\ne\ell$, then
	\[
		\frac{\th_r-\th_s}{\th_k-\th_\ell} \in \rats.
	\]
	\item
	The vertices $a$ and $b$ are strongly cospectral.
\end{enumerate}
We refer to (a), (b) respectively as symmetry and periodicity, (c) is related to what is
known as the ratio condition. 

One goal of this paper is to show that all these properties are consequences of the fact
that the density matrices $D_a$ and $D_b$ are real. Thus if $D_1$ and $D_2$ are real
density matrices and $D_1(t)=D_2$, then strong analogs of the above four properties hold.
(In fact these properties will be easy corollaries of our more general results.)

We will also see that interesting things happen if we assume that the entries of $D_1$
and $D_2$ are algebraic numbers.

For background on state transfer and continuous quantum walks, 
see e.g. \cite{st-transfer, kay2010perfect}.

\section{Real State Transfer}

A state is \textsl{real} if its density matrix is real. We recall that we have
perfect state transfer at time $t$ from a density matrix $P$ to a density matrix $Q$
if $P(t)=Q$. We say that we have \textsl{pretty good state transfer}
from a state $P$ to $Q$ if for each positive real number $\eps$ there
is a time $t_{\epsilon}$ such that $\norm{(U(t)PU(t)-Q)}<\epsilon$.

We define the \textsl{eigenvalue support of a density matrix} $P$ to be the set
of pairs $(\theta_r,\theta_s)$ such that $E_rPE_s\ne0$. (This definition extends the one 
given in the introduction to density matrices not of the form $D_a$.)

\begin{lemma}
	If $P$ is a density matrix and $E_r, E_s$ are spectral idempotents of $A$
	such that $E_rPE_s\ne0$ then neither $E_rPE_r$ nor $E_sPE_s$ is zero.
\end{lemma}

\proof
Since $P$ is positive semidefinite, there is a unique positive semidefinite matrix
$Q$ such that $P=Q^2$. Hence if $E_rPE_s\ne0$, then $E_rQ\ne0$ and $E_sQ\ne0$.
Hence $E_rPE_r = E_rQ^2E_r\ne0$ and similarly $E_sPE_s\ne0$.\qed

If $E_rPE_s=0$ whenever $r\ne s$, then
\[
	P = \sum_r E_rPE_r
\]
and $PA=AP$; thus $U(t)PU(-t)=P$ for all $t$.

We say that the eigenvalue support of $P$ satisfies the ratio condition if, 
for any two pairs of distinct eigenvalues $(\th_r,\th_s)$ and $(\th_k,\th_ell)$ 
in the eigenvalue support of $P$ with $\th_k\ne\th_\ell$, we have
\[
	\frac{\th_r-\th_s}{\th_k-\th_\ell} \in \rats.
\]
Note that if $P$ is pure, that is, $P=xx^*$ for some unit vector $x$, then
$E_rPE_s=0$ if and only if $E_rx$ or $E_sx$ is zero. (In this case we could define
the eigenvalue support to be the the set of eigenvalues $\th_r$ such that
$E_rPE_r\ne0$, which what we do elsewhere.)

\begin{theorem}\label{thm:real-return}
	Let $U(t)$ be the transition matrix corresponding to a graph $X$.
	Let $=\sum_r e^{it\th_r}E_r$ be the spectral decomposition of $A$.
	If $P$ is a real density matrix, there is a positive time $t$ such that
	$U(t)PU(-t)$ is real if and only if the eigenvalue support of $P$ satisfies the
	ratio condition.
\end{theorem}

\proof
We have
\[
	P^{(t)} = \sum_{r,s}e^{it(\th_r-\th_s)} E_rPE_s
\]
and therefore the imaginary part of $P^{(t)}$ is
\[
	\sum_{r,s} \sin(t(\th_r-\th_s)) E_rPE_s.
\]
The non-zero matrices $E_rPE_s$ are linearly independent, and therefore the imaginary part
of $P^{(t)}$ is zero if and only if $\sin(t(\theta_r-\theta_s))=0$ whenever $E_rPE_s\ne0$.
Hence if $(\theta_r,\theta_s)$ lies in the eigenvalue support of $P$ then there is an
integer $m_{r,s}$ such that $t(\theta_r-\theta_s)=m_{r,s}\pi$.\qed

If $E_rPE_s=0$ whenever $r\ne s$, then
\[
	P = \sum_r E_rPE_r
\]
and $PA=AP$; thus $U(t)PU(-t)=P$ for all $t$.

\begin{lemma}\label{lem:period}
    Let $P$ be a real density matrix. If $P(t)$ is real, then $P(2t)=P$
	and $U(2t)$ commutes with $P$ and $Q$.
\end{lemma}

\proof
If $U(t)PU(-t)=Q$ where $Q$ is real, then taking complex conjugates yields 
\[
	U(-t)PU(t)=Q
\]
and consequently $P=U(t)QU(-t)$. It follows at one that $U(2t)$ commutes
with $P$.\qed

\begin{lemma}
    If $A$ is real and symmetric and we have state transfer at $t$ between
    real density matrices $P$ and $Q$, then 
    \begin{enumerate}[(a)]
        \item 
        $E_rPE_r = E_rQE_r$.
        \item
        If $t(\theta_r-\theta_s)$ is not a multiple of $\pi$ then $E_rPE_s=E_rQE_s=0$.
        \item
        Otherwise $E_rPE_s = \pm E_rQE_s$.
    \end{enumerate}
\end{lemma}

\proof
If $A$ is real and symmetric then the idempotents $E_r$ are real and symmetric.
Assume
\[
	Q = U(\tau)P U(-\tau) = \sum_{r,s} e^{i\tau(\theta_r-\theta_s)}E_r P E_s.
\]
If we multiply this expression on the left by $E_r$ and on the right by $E_s$, then
\[
    E_r Q E_s = e^{i\tau(\theta_r-\theta_s)}E_r P E_s
\]
and, since all matrices here are real, $e^{it (\theta_r-\theta_s)}$ must be real.

Since $\sum_{r,s}E_rPE_s=P$, we see that if $A$ is real and there is pst from 
$P$ to $Q$,
there are signs $\eps_{r,s}=\pm1$ such that
\[
    Q = \sum_{r,s} \eps_{r,s} E_rPE_s.\qed
\]

Consider the rank-one case. If $P=uu^T$ then $E_ruu^TE_s=0$ if and only if $E_ru=0$
or $E_su=0$. Hence the constraint in (b) gives a constraint on the eigenvalue support
of $u$. (In fact this is the usual ratio condition, so (b) generalizes this.)

\section{Pretty Good State Transfer}
\label{sec:pgst}

We have \textsl{pretty good state transfer} from a state $P$ to a state $Q$ if,
for each $\eps>0$, there is a time $t$ such that
\[
	\norm{U(t)PU(-t)-Q}<\eps.
\]
Since the complex conjugate of 
\[
	U(t)PU(-t)-Q
\]
is
\[
    U(-t)PU(t)-Q = U(-t)\,(P-U(t)QU(-t))\,U(t)
\]
and since $U(t)$ is unitary,
\[
    \norm{U(t)PU(-t)-Q} = \norm{P-U(t)QU(-t)}.
\]
Hence if we have pretty good state transfer from $P$ to $Q$, then we have pretty 
good state transfer from $Q$ to $P$.

\begin{lemma}
    Suppose $A$ is real and we have pretty good state transfer from $P$ to $Q$. Then 
    $E_rPE_s = \pm E_rQE_s$ (for all $r$ and $s$) and $E_rPE_r=E_rQE_r$.
\end{lemma}

\proof
By assumption there is an increasing sequence of times $(t_k)_{k\ge0}$ such that 
\[
    U(t_k)PU(-t_k) \to Q
\]
as $t_k\to\infty$. Hence
\[
    e^{i(\theta_r-\theta_s)t_k} E_rPE_s \to E_rQE_s
\]
as $t_k\to\infty$. Since $E_rPE_s$ and $E_rQE_s$ are real, the result follows.\qed

\begin{corollary}
    If $A$ is real and $P$ is real, there are only finitely many real density 
	matrices $Q$ for which there is pretty good state transfer from $P$ to $Q$.
\end{corollary}

\proof
Since $Q=\sum_{r,s}E_rQE_s$ it follows that $Q=\sum_{r,s}\epsilon_{r,s}E_rPE_s$,
where $\epsilon_{r,s}=\pm1$.\qed

Pretty good state transfer is treated in some detail in \cite{path-pgst}.

\section{Algebras}

Because they are trace-orthogonal, the non-zero matrices $E_rPE_s$ are linearly independent. 
The ``off-diagonal'' terms $E_rPE_s$ are nilpotent, indeed the matrices
\[
	E_rPE_s,\quad (r<s)
\]
generate a nilpotent algebra where the product of any two elements is zero. 
The ``diagonal'' terms $E_rPE_r$ generate a commutative semisimple algebra; their sum
is the orthogonal projection of $P$ onto the commutant of $A$. (See \cite[Section~5]{cgjs-strong} 
for further details).

Since $U(t)$ is a linear combination of the spectral idempotents of $A$, it is a polynomial
in $A$ and therefore, for each $t$ we that $P(t)\in\alg{A,P}$. In consequence
\[
	\alg{P(t),A} = \alg{P,A}
\]
for all $t$.

\begin{lemma}
    If we have pretty good state transfer from $P$ to $Q$, then 
	$\alg{A,P}=\alg{A,Q}$.\qed
\end{lemma}

\proof
The algebra $\alg{A,P}$ is closed and as $Q$ is a limit of a sequence of matrices
in it, it follows that $Q\in\alg{A,P}$. If we have pretty good state transfer from
$P$ to $Q$, then as we noted st the beginning of Section~\ref{sec:pgst}, 
there is pretty good state transfer from $Q$ to $P$ and so $\alg{A,P}=\alg{A,Q}$.\qed

If $\alg{A,P}$ is the full matrix algebra, we say that $P$ is \textsl{controllable}.
If $P$ is real and $P(t)$ is real, then $U(2t)$ commutes with
$A$ and $P$. If $P$ is also controllable it follows that $U(2t)$ must be a scalar matrix,
say $U(2t)=\zeta I$. Since $\det(U(t))=1$, we have $\zeta^{|V(X)|}=1$ and therfore
$\zeta$ is a root of unity.

\section{Timing}

We investigate the times at which perfect state transfer can occur.

\begin{lemma}
	Let $P$ be a density matrix and let $S$ be given by
	\[
		S := \{\sigma: U(\sigma)PU(-\sigma)=P\}.
	\]
	Then there are three possibilities:
	\begin{enumerate}[(a)]
		\item 
		$S=\emptyset$.
		\item
		There is a positive real number $\tau$ and $S$ consists of all integer multiples
		of $\tau$.
		\item
		$S$ is a dense subset of $\re$ and $U(t)PU(-t)=P$ for all $t$.
	\end{enumerate}
\end{lemma}

\proof
Suppose $S\ne\emptyset$.
Then $S$ is an additive subgroup of $\re$ and there are two cases.
First, $S$ is discrete and consists of all integer multiples of its least 
positive element. 
Second, $S$ is dense in $\re$ and there is sequence of positive elements
$(\sigma_i)_{i\ge0}$ with limit 0. Since for small values of $t$ we have
\[
	U(t) \approx I+itA
\]
it follows that $AP=PA$ and $U(t)PU(-t)=P$ for all $t$.\qed

If $U(t)PU(-t)=P$ when $t=\tau>0$, but not for all $t$, we say that $P$ is 
\textsl{periodic} with period $\tau$. If a density matrix is periodic, 
it has a well defined 
minimum period. If there is perfect state transfer from $P$ to $Q$, then $P$
and $P$ are both periodic with the same minimum period.

\begin{lemma}
    Suppose $P$ and $Q$ are distinct real density matrices.  If there is 
	perfect state transfer from $P$ to $Q$, then $P$ is periodic and 
	if the minimum period of $P$ is $\sigma$, then pst occurs at time $\sigma/2$.
\end{lemma}

\proof
Suppose we have pst from $P$ to $Q$ and define
\[
    T := \{t: U(t)PU(-t)=Q\}
\]
Assume that the minimum period of $P$ is $\sigma$. If $t\in T$ then $P$ is periodic
with period $2t$ and so $T$ is a coset of a discrete subgroup of $\re$. Also $T=-T$.
Let $\tau$ be the least positive element of $T$. Then $2\tau\ge \sigma$ and thus
\[
    \tau \ge \frac12 \sigma.
\]
If $\tau\ge \sigma$ then $\tau-\sigma\in T$ and since $\tau$ is not a period,
$\tau<\sigma$. As $\sigma$ must divide $2\tau$, it follows that $\tau=\sigma/2$.\qed

\begin{corollary}
    For any real density matrix $P$, there is at most one real density 
	matrix $Q$ such that there is perfect state transfer from $P$ to $Q$.\qed
\end{corollary}

\begin{lemma}
    Suppose we have pst from $P$ to $Q$ at time $t$ and that $\seq \theta1m$ are the
    distinct eigenvalues of $A$ in nonincreasing order. If $\tr(PQ)=0$ then
    \[
        t \ge \frac{\pi}{\theta_1-\theta_m}.
    \]
\end{lemma}

Suppose $U(t)PU(-t)=Q$ and $\tr(PQ)=0$. We have
\[
    U(t)PU(-t) =\sum_r e^{i \theta_r t} U(t)E_rU(-t)
\]
and therefore
\begin{align*}
    \tr(PQ) &= \sum_r e^{i \theta_r t} \tr(PU(t)E_rU(-t))\\
        &= \sum_r e^{i \theta_r t} \tr(U(-t)PU(t)E_r)\\
        &= \sum_r e^{i \theta_r t} \tr(QE_r)
\end{align*}
Since $Q$ and the idempotents $E_r$ are positive semidefinite, $\tr(QE_r)\ge0$.
Also
\[
    \sum_r \tr(QE_r) = \tr(Q) = 1.
\]
Hence if $\tr(PQ)=0$ then we see that 0 is a convex combination of the eigenvalues
$e^{i \theta_rt}$ of $U(t)$. If $A$ has $m$ distinct eigenvalues
\[
    \theta_1\ge\cdots\ge \theta_m.
\]
this implies that $e^{it \theta_1},\ldots,e^{it \theta_m}$ cannot be contained in an arc 
on the unit circle in the complex plane of length less than $\pi$, and therefore
$t(\theta_1-\theta_m)\ge\pi$.\qed

Note that in this lemma we do not need $P$ and $Q$ to be real.

An algebraic integer is \textsl{totally real} if all its algebraic conjugates are
real, equivalently, all zeros of its minimal polynomial are real.

\begin{theorem}
	Let $P$ be a rational state with eigenvalue support $S$. If $S$
	satisfies the ratio condition, then there is a square-free integer $\De$
	such that if $E_rPE_s\ne0$, then $\th_r-\th_s$ is an integer multiple of $\sqrt{\De}$.
\end{theorem}

\proof
Let $\flde$ denote the extension field of $\rats$ generated by the elements of $S$.
Since $P$ is rational, it follows that if $E_rPE_s\ne0$ and $\ga\in\Ga$, then
\[
	E_r^\ga P E_r^\ga \ne 0.
\]
We also note that, as $A$ is an integer matrix, the spectral
idempotents $E_r$ are algebraic and therefore $E_r^\ga$ is a spectral idempotent of $A$.

The product 
\[
	\prod_{(\th_r,\th_s)\in S} \th_r-\th_s
\]
is fixed by $\Ga$ and is consequently an integer. Given the ratio condition, we see that
if $(\th_k,\th_\ell)\in S$, then
\[
	\prod_{(\th_r,\th_s)\in S} \frac{\th_k-\th_\ell}{\th_r-\th_s} \in \rats
\]
and therefore
\[
	(\th_k-\th_\ell)^{|S|} \in \rats.
\]
As the eigenvalues of $A$ are algebraic integers, this implies that $(\th_k-\th_\ell)^{|S|}$
is an integer. The eigenvalues of $A$ are totally real, but if the polynomial $t^s-1$
has a real root, then if must be $\pm1$ and in this case $s$ must even. We conclude that
$(\th_k-\th_\ell)^2$ is an integer.

Suppose there are integers $m_{k,\ell}$ and $m_{r,s}$ and square-free integers $b$ and $c$
such that
\[
	\th_k-\th_\ell = m_{k,\ell}\sqrt{b},\qquad \th_r-\th_s = m_{r,s}\sqrt{c}.
\]
If 
\[
	\frac{\sqrt{b}}{\sqrt{c}} = \frac{\th_k-\th_\ell}{\th_r-\th_s} \in \rats,
\]
tne $b=c$.\qed

\begin{corollary}
	If $P$ is a periodic rational state, then the period of $P$ is at most $2\pi$.
\end{corollary}

\proof
If $t=2\pi/\sqrt{\De}$ then $t(\th_r-\th_s)$ is an even multiple of $\pi$ for each
pair $(\th_r,\th_s)$ in the eigenvalue support of $P$. Consequently
\[
	P^{(t)} = \sum_{r,s} e^{it(\th_r-\th_s)}E_rPE_s = \sum_{r,s} E_rPE_s = P.\qed
\]

\section{Algebraic States}

We say that a state with density matrix $D$ is \textsl{algebraic} if the
entries of $D$ are algebraic numbers. Clearly the states $D_a$ are algebraic.

Suppose $D$ is a pure state. Then $D(t)$ is pure for all $t$. If $D=zz^*$,
then $D(t)=ww^*$, where $w=U(t)z$. We say that a matrix or vector is
\textsl{flat} if all its entries have the same absolute value. We see
that a vector $w$ is flat if and only the diagonal entries of $ww^*$
are all equal. (Note that these entries are non-negative and real.)

We say that a quantum walk has \textsl{uniform mixing relative to
a pure state} $D$ if there is a time $t$ such that
\[
	D(t)\circ I = \frac1n I.
\]
We refer to uniform mixing relative to the vertex state $D_a$ as \textsl{local uniform mixing}.
The walk has \textsl{uniform mixing} if, it at some $t$, it admits uniform mixing relative
to each vertex. In many of the cases where uniform mixing is known to occur,
the underlying graph is vertex transitive, and then uniform mixing occurs
if and only if uniform mixing relative to a vertex occurs. The only
examples we know of graphs that are not regular and that do admit uniform mixing are
the complete bipartite graph $K_{1,3}$ and its Cartesian powers (an observation due to H. Zhan).
If $n\ge2$, the stars $K_{1,n}$ admit uniform mixing relative to the
vertex of degree $n$.

Carlson et al.~\cite{carlson-unimix} observed that we have uniform mixing relative to the vertex 
of degree $n$ in the star $K_{1,n}$.

We say that the continuous quantum walk on $X$ is \textsl{periodic at $a$}
if there is a time $\tau$ such that $D_a(\tau)=D_a$.

\begin{theorem}
	\label{thm:alg-ratio}
	Let $A$ be a Hermitian matrix with algebraic entries and let $U(t)=\exp(itA)$
	If the density $D$ is algebraic and, for some $t$, the density $D(t)$
	is algebraic, then the ratio condition holds on the eigenvalue support of $D$.
\end{theorem}

\proof
Since the entries of $A$ are algebraic, its eigenvalues are algebraic and therefore
the spectral idempotents are algebraic.

We have
\[
	D(t) = \sum_{r,s}e^{it(\th_r-\th_s)}E_rDE_s.
\]
The matrices $E_rDE_s$ are pairwise orthogonal, and so, for all $r$ and $s$,
\[
	\ip{D(t)}{E_rDE_s} = e^{it(\th_r-\th_s)} \ip{E_rDE_s}{E_rDE_s}.
\]
The entries of the spectral idempotents are algebraic, and if the entries
of $D$ and $D(t)$ are algebraic, then the values of the two inner products
in the above identity are algebraic numbers.

It follows that $e^{it(\th_r-\th_s)}$ must be algebraic, for all $r$ and $s$.
Now if $k\ne\ell$, then
\[
	e^{it(\th_r-\th_s)} 
		= \left(e^{it(\th_k-\th_\ell)}\right)^{\frac{\th_r-\th_s}{\th_k-\th_\ell}}.
\]
The Gelfond-Schneider theorem tells us that if $\al$ and $\be$ are algebraic numbers
and $\al\ne0,1$ and $\be$ is irrational, then $\al^\be$ is transcendental, whence
we deduce that if $D$ and $D(t)$ are algebraic, then if $k\ne\ell$ and neither $E_rDE_s$
nor $E_kDE_\ell$ is zero, the ratio
\[
	\frac{\th_r-\th_s}{\th_k-\th_\ell}
\]
is rational.\qed

The Gelfond-Schneider theorem was first used as above in \cite{adamczak2007}; the technique 
is due to Jennifer Lin. One source for the Gelfond-Schneider theorem is Burger and 
Tubbs \cite{burger2004}.

\section{Oriented Graphs}

We study state transfer on oriented graphs. In this section we consider the graph theory
and linear algebra, in the next we turn to the continuous walks.

An oriented graph is a structure consisting of vertices and arcs, where an arc is an ordered
pair of vertices, and any two vertices lie in at most one arc. We can construct (a large number of)
oriented graphs by choosing a graph and assigning a direction to each edge. We use $\arcs(X)$
to denote the set of arcs of $X$. The adjacency matrix $S$ of $X$ id the matrix with rows and 
columns indexed by $V(X)$, where
\[
	S_{u,v} = \begin{cases}
		1,& uv \in \arcs(X);\\
		-1,& vu\in\arcs(X);\\
		0,& \mathrm{otherwise}.
		\end{cases}
\]
Thus $S$ is a skew symmetric matrix. We define the \textsl{degree} of a vertex $v$ in $X$
to be the number of arcs that use $v$. As we defined them, each oriented graph has an underlying
graph whose adjacency matrix is $S\circ S$. The total valency of a vertex in $X$ is its valency
in the undirected graph that underlies $X$.

The matrix $iS$ is Hermitian, with all eigenvalues real, and therefore the eigenvalues of $S$ 
are purely imaginary. They are symmetric about the real axis of the
complex plane, so we will assume that the $r$-th eigenvalue is $i\la_r$, where $\la_r$ is real. 
We can then write the spectral decomposition of $S$ as
\[
	S = \sum i\la_r F_r
\]
where the idempotents $F_r$ are Hermitian. Further, the idempotent associated to the eigenvalue
$-i\la_r$ is $\comp{F}_r$.

\begin{lemma}
	Let $X$ be an oriented graph with maximum total valency $\De$. If $\la$ is an eigenvalue of
	$X$, then $|\la|\le\De$.
\end{lemma}

\proof
If $z\ne0$ and $Sz=\la z$, then
\[
	\la z_k = \sum_{\ell} S_{k,\ell}z_\ell
\]
and so by the triangle inequality,
\[
	\abs{\la} \abs{z_k} \le \sum_{\ell: S_{k,\ell}\ne0} |z_\ell|.
\]
By choosing $k$ so that $\abs{z_k}$ is maximal, we obtain the stated bound.\qed

If $Y$ is a bipartite graph and
\[
	A(Y) = \pmat{0&B\\ B^T&0},
\]
then
\[
	S = \pmat{0& -B\\ B^T&0}
\]
is skew symmetric. As
\[
	\pmat{-iI&0\\ 0&I}\pmat{0& -B\\ B^T&0}\pmat{iI&0\\ 0&I} = i\pmat{0&B\\ B^T&0},
\]
each eigenvalue of $S$ is equal to $i\th$, for some eigenvalue $\th$ of $Y$. The oriented
graph with adjacency matrix $S$ will be called the \textsl{natural orientation} of $Y$.
In a sense, the spectral
theory of bipartite graphs is the intersection of the spectral theory of graphs with the
spectral theory of oriented graphs.

Finally, since $iS$ is Hermitian, the eigenvalues of a principal submatrix interlace the eigenvalues
of $S$, and therefore the eigenvalues of an induced subgraph of an oriented graph $X$ interlace
the eigenvalues of $X$.

\section{Quantum Walks on Oriented Graphs}

Continuous quantum walks on oriented graphs were first studied in \cite{cameron2014universal, connelly2017universality}.

In the introduction we defined the transition matrix $U(t)$ as $\exp(itA)$, where $A$ was 
the adjacency matrix of a graph. Thus $A$ was symmetric and real. 
However all that is needed is that $A$ 
should be Hermitian and hence, if $S$ is the adjacency matrix of an oriented graph,
we may define a transition matrix
\[
	U(t) = \exp(it(-iS)) = \exp(tS).
\]
We note that $U(t)$ is then real and orthogonal for all $t$. We have the spectral decomposition
\[
	U(t) = \sum_r e^{it\la_r}F_r
\]
where $\la_r\in\re$ and $F_r$ is Hermitian. 

As we noted in the previous section, if $F$ is the spectral idempotent associated to the eigenvalue
$i\la$, then the eigenvalue associated to $-i\la$ is $\comp{F}$. Hence $F_re_a=0$
if and only $\comp{F}e_a=0$. Consequently the eigenvalue support of a vertex is symmetric about 
the real axis of the complex plane and therefore the ratio condition on the eigenvalue support of
a vertex is equivalent to the condition that $\la/\mu\in\rats$ for all choices of $\la$ and $\mu$.

If $S$ arises as the adjacency matrix of the natural orientation of a bipartite graph $Y$,
with adjacency matrix $A$, then
\[
	\pmat{-iI&0\\ 0&I} \exp(tS) \pmat{iI&0\\ 0&I} = \exp(itA).
\]
Accordingly we have perfect state transfer from $a$ to $b$ relative to $\exp(itA)$ if and only if 
it occurs from $a$ to $b$ relative to $\exp(tS)$; similarly we have local uniform mixing at
$a$ in the graph if and only if we have it in the oriented graph.

\begin{theorem}
	If there is perfect state transfer on an oriented graph from a vertex $a$, then the
	eigenvalue support of $a$ satisfies the ratio condition.
\end{theorem}

\proof
We simply note that $D_a$ and $D_b$ are algebraic, whence Theorem~\ref{thm:alg-ratio}
implies the conclusion.\qed

\begin{theorem}
	If there is local uniform mixing at a vertex $a$ in an oriented graph, then the
	eigenvalue support of $a$ satisfies the ratio condition.
\end{theorem}

\proof
Assume $n=\abs{V(X)}$.
If there is local uniform mixing at $a$ at tine $t$, then $U(t)e_a$ is flat. As $U(t)$ is real,
this implies that the entries of $U(t)e_a$ are all equal to $\pm n^{1/2}$ and hence they are
algebraic. Now apply Theorem~\ref{thm:alg-ratio}.\qed

We say that an oriented graph is \textsl{connected} if its underlying undirected graph is connected.

\begin{lemma}
	Let $X$ be an oriented graph with adjacency matrix $S$ and suppose $a\in V(X)$. 
	Then $\exp(tS)e_a$ is periodic if and only if the ratio condition holds on the
	eigenvalue support of $a$.
\end{lemma}

\proof
Suppose $U(t)e_a=e_a$. Then
\[
	e_a = U(t)e_a = \sum_e e^{it\la_r} F_re_a
\]
and since we also have
\[
	e_a = \sum_r F_re_a,
\]
it follows that $e^{it\la_r}=1$ for all $r$ such that $i\la_r\in\esupp(a)$. Consequently
$t\la_r$ is an integer multiple of $2\pi$ for each $r$, and therefore the ratio of any two
elements of $\esupp(a)$ is rational.

Now assume conversely that the ratio condition holds on the eigenvalue support of $a$, and
set $m=\abs{\esupp(a)}$. Let $\Ga$ denote the Galois group of the extension field 
of $\rats$ generated by $\seq\la1s$. Then $\esupp(a)$ is closed under $\Ga$. Therefore
\[
	\prod_{r=1}^m \la_r
\]
is fixed by each element of $\Ga$, and is thus an integer. Hence $\prod_r\la_r\in\ints$. Now
\[
	\prod_{r=1}^m \frac{\la_s}{\la_r}
\]
is rational and accordingly $\la_s^m$ is rational. Since it also an algebraic integer, 
it must be an integer. Because it is an eigenvalue of the Hermitian matrix $iS$, all algebraic conjugates
of $\la_r$ are real, and therefore we must have $\la_r^2\in\ints$.Therefore for each $s$ we have
$\la_s = a_s\sqrt{b_s}$ where $a_s,b_s\in\ints$ and $b_s$ is square free. Since $\la_s/\la_r$
is rational it follows that $b_s$ is independent of $s$. Therefore $\esupp(a)$ consists of
integer multiples of $\sqrt{b}$, for some square-free integer $b$. This shows that $a$ is periodic,
with period $2\pi/\sqrt{b}$.\qed

\begin{theorem}
	There are only finitely many connected bipartite graphs with maximum valency at most $k$
	which contain a periodic vertex.
\end{theorem}

\proof
Let $c$ be the eccentricity of $a$. The vectors 
\[
	(A+I)^re_a, \qquad (r=0,\ldots,c) 
\]
are linearly independent, because their supports form
a strictly increasing sequence. These vectors lie in the span of the non-zero vectors $F_re_a$,
whence $c+1$ is bounded above by the size of the eigenvalue support of $a$. Since any two elements
of $\esupp(a)$ differ by a least one, and since the largest eigenvalue of $S$ is at most $k$,
we have that $\abs{\esupp(a)}\le2\De+1$. Hence the eccentricity of $a$ is bounded by a function
of $\De$, and therefore $|V(X)|$ is bounded.\qed

\begin{corollary}
	There are only finitely many connected bipartite graphs with maximum valency at most $k$
	which admit local uniform mixing.\qed
\end{corollary}

The theorem also implies that there are only finitely many connected bipartite graphs with 
maximum valency at most $k$ on which perfect state transfer occurs, but this holds more
generally for all graphs, bipartite or not. See \cite[Corollary 6.2]{when-pst}.

\section{Prospects, Problems}

We have established a surprising connection between behaviour
of continuous quantum walks at certain times and the field of definition of
the associated density matrix. An obvious problem is to find more examples
of such behaviour.

Secondly, most work on continuous quantum walk assumes that the initial state
is of the form $e_re_r^T$. Our results indicate that it might be fruitful
to consider more general initial states. Our personal feeling is that
pure states will be most interesting, because their eigenvalue support tends
to be smaller.


\end{document}